\documentclass[12pt,A4paper]{article}

\usepackage[left=1in, right=1in, top=1in, bottom=30mm]{geometry}

\usepackage{epsfig}
\usepackage{epstopdf}
\usepackage{soul}
\usepackage{xfrac}
\usepackage{amsmath}
\usepackage{enumerate}

\usepackage{subcaption}

\usepackage[percent]{overpic}

\usepackage{color}
\usepackage{amsthm,amsmath,amssymb,mathrsfs}
\usepackage{booktabs}
\usepackage{mathpazo}
\usepackage{microtype}
\usepackage{overpic}
\usepackage{bm}
\usepackage{sectsty}
\usepackage[
	pdftitle={PDFTitle},
	pdfauthor={Hugo Parlier},
	ocgcolorlinks,
	linkcolor=linkred,
	citecolor=linkred,
	urlcolor=linkblue]
{hyperref}

\usepackage{tikz}
\usetikzlibrary{calc,decorations.pathreplacing}

\definecolor{linkred}{RGB}{112, 11, 75} 
\definecolor{linkblue}{RGB}{16, 78, 139}

\usepackage[hang,flushmargin]{footmisc}

\usepackage{titlesec}
	\titlespacing{\section}{0pt}{12pt}{0pt}
	\titlespacing{\subsection}{0pt}{6pt}{0pt}
	
\titlelabel{\thetitle.\quad}

\makeatother 

\theoremstyle{plain}
\newtheorem{theorem}{Theorem}[section]

\newtheorem{lemma}[theorem]{Lemma}

\theoremstyle{definition}

\newtheorem{remark}[theorem]{Remark}

\newcommand{\arc}{{\mathcal A}}

\newcommand{\M}{{\mathcal M}}

\newcommand{\area}{{\rm area}}
\newcommand{\arcsinh}{{\,\rm arcsinh}}

\sectionfont{\large \bfseries}
\subsectionfont{\normalsize}

\long\def\symbolfootnote[#1]#2{\begingroup%
\def\thefootnote{\fnsymbol{footnote}}\footnote[#1]{#2}\endgroup}

\def\blfootnote{\xdef\@thefnmark{}\@footnotetext}

\usepackage{mathtools}

\usepackage[sort,nocompress]{cite}

\begin{document}

{\Large \bfseries The full length of short pants}

{\large Hugo Parlier\symbolfootnote[1]{\small Supported by ANR-SNF Grant number 200021E\_238147
(SUGAR).
}
}

\vspace{0.5cm}
{\bf Abstract.}
Bers proved that any genus $g$ hyperbolic surface admits a pants decomposition of length bounded above by a function of its genus. The goal here is to provide an explicit upper bound for the sums of lengths which is on the order of $(g \log(g)) ^{3/2}$.
\vspace{0.5cm}

\section{Introduction}\label{sec:intro}
The fact that any closed hyperbolic surface admits a pants decomposition of length bounded solely by its topology is particularly useful in the study of the deformation spaces of surfaces. Since this was first observed by Bers \cite{Bers1,Bers2}, the optimal length bounds are generally called Bers' constants.

There are in fact multiple ways of measuring the lengths of pants decompositions, by considering for instance the maximal length, or the sum of the lengths, and hence there are different ways of defining Bers' constants. After some initial quantifications in the 1970s, in the 1980s Buser famously studied the maximal length and showed that there always exist a pants decomposition where all lengths are of length at most roughly $g$ where $g$ is the genus of the surface (see \cite{BuserPants,BuserSeppala, BuserBook}). Since there are $3g-3$ curves in a pants decomposition, this immediately implies that the sums of lengths are bounded above by  $\approx g^2$. The constants in front of the genus have been improved (see \cite{ParlierShort,ParlierShorter}) but the best upper bounds, up to a logarithmic factor, have remained on these orders since the 1980s. In parallel, Buser also produced families of surfaces where every pants decomposition necessarily has a curve of length approximately $\sqrt{g}$ (the so-called hairy torus examples), but it wasn't until the Buser-Sarnak examples of surfaces \cite{BuserSarnak94} where every curve is of length at least $\approx \log(g)$ that a supra-linear example for the sums of lengths was found. While, to this date, those examples remain the best explicit lower bound examples for the sum, random surfaces, with respect to the Weil-Petersson volume form, were shown to have all pants decompositions of total length at least $\approx g^{7/6 -\varepsilon}$ for any $\varepsilon$ and large enough genus \cite{GPY}. 

Buser also conjectured that the optimal bound for the maximal length should be on the order of $\sqrt{g}$. In a similar direction, Balacheff and Sabourau proved a related diastolic inequality \cite{BalacheffSabourau}, which can be used to prove another linear upper bound, see for instance \cite{CavendishParlier}, and also inspired Guth to propose a variation on Buser's problem as a question about the length of fibers of maps from surfaces to graphs \cite{GuthICM2010}. In any event, if Buser's conjecture is correct, that would imply an upper bound on the sums of the lengths on the order of $g^{3/2}$, but, as stated before, the known upper bounds for the max remain stuck at linear in the genus.

If, instead of increasing topology by adding genus, one increases topology by leaving the closed surface setup and adding punctures, then the situation is quite different. For punctured spheres, with $n$ punctures, the corresponding question for the max was resolved \cite{BalacheffParlier}, showing there are pants decompositions with all curves of length at most $\sqrt{n}$, hence solving Buser's conjecture in that setup. For total length however, an upper bound was shown, together with Balacheff and Sabourau, which is on the order of $n \log(n)$, which beats both the direct $n^{3/2}$ bound, and shows that for sums of lengths, punctures behave quite differently from the genus bounds \cite{BPS2012}. 

The main goal of this paper is to show the following upper bound, which provides the first notable order of growth improvement for closed surfaces since the 1980s. 

\begin{theorem}\label{thm:pantsintro}
There exists a universal constant $A$ such that every closed hyperbolic surface $X$ of genus $g\geq2$ admits a pants decomposition of total length at most
$
A\,g^{3/2}(\log g)^{3/2}.
$
\end{theorem}

The constant $A$ can be taken to be $80000$ (see Theorem \ref{thm:pants}), but no specific effort has been made to optimize it. The main technical step to obtain this result concerns the size of a cut system, that is a collection of $g$ disjoint simple closed curves such that the surface obtained by cutting along them is planar (and hence a $2g$-holed sphere). 

\begin{theorem}\label{thm:cutsystemintro}
There exists a universal constant $B$ such that every closed hyperbolic surface $X$ of genus $g\geq 2$ admits a cut system of total length at most
$B\,g^{3/2}\sqrt{\log g}$.
\end{theorem}
Again, an explicit bound of $3000$ for the constant $B$ is provided in Theorem \ref{thm:cutsystem}, but obviously there has been little attempt at making this optimal. 

A word about proofs: using \cite{BPS2012}, a method to get from Theorem \ref{thm:cutsystemintro} to Theorem \ref{thm:pants} is already quite clear, so the real novelty is in the cut system result. Now the idea to get the cut system result is also inspired by this same result. The first step is to create a combinatorial model for the surface, namely a cell decomposition of the surface with bounded edge length, making amends for the thin parts of the surface. Here is where the strategy becomes more surprising: instead of looking for short non-trivial cycles in the embedded graph, by looking at the neighborhood of a spanning tree of the graph, we can construct a curve that is uniformly dense on the thick part of the surface. This curve has multiple obvious defaults: it is long and, if the surface is uniformly thick, it is trivial. And yet it has a very strong redeeming feature: all remaining edges of the graph can be used to create curves with it. In the case of punctured surfaces \cite{BPS2012}, the corresponding "filling" curve was non-trivial and we used this curve to surround a full pants decomposition. Here the topology makes it difficult to find disjoint curves, but what you can do is find a batch of roughly $\sqrt{g \log(g)}$ curves that can be used to remove a planar subsurface of a certain size in a way that can be repeated. 

Finally, we mention that these methods are transposable to the obvious generalizations, including Riemannian metrics, non-orientable surfaces, or combinatorial settings such as in \cite{CHdM}. The question that remains, for sums of lengths, is the correct exponent which is now known to lie somewhere between $\frac{7}{6}$ and $\frac{3}{2}$. This is not unlike the problem for the maximal systole growth problem, where the remaining question is about the term in front of the logarithm which has been known to be somewhere between $\frac{4}{3}$ and $2$ since the work of Buser and Sarnak \cite{BuserSarnak94}.

\noindent {\bf Acknowledgements.}

This work benefited greatly from many conversations about pants decompositions over the years. I thank my collaborators on these topics Florent Balacheff, Peter Buser, Larry Guth, Stéphane Sabourau, and Robert Young in particular.

\noindent {\bf Statement on the use of AI.}

While preparing the manuscript, certain standard AI tools were used for formatting, bibliographical searches and checking certain computations. In addition, one conversation with an LLM led to Lemma \ref{lem:circle} which turns out to be crucial to the proof. I wanted to use a type of Ramsey theorem to construct many curves simultaneously (or a large planar surface in one shot) but the traditional Ramsey theorem only allows for $\approx \log(g)$ curves. The LLM correctly pointed out that the particular structure of circle graphs allows for many more.

\section{Preliminaries}\label{sec:prelim}

We begin with the definitions and notation that will be used throughout.

\subsection{Definitions and notation}
Surfaces are orientable and closed. Moduli space $\M_g$ is the space of hyperbolic surfaces of genus $g\geq 2$, up to isometry. Throughout, $X\in \M_g$ will be a closed surface, and by an Euler characteristic argument is of area $4\pi (g-1)$. 

Closed curves on $X$ are considered to be non-oriented and are called simple if there are no double points. Closed curves on $X$ are either trivial, that is homotopic to a point, or homotopic to a unique closed geodesic, which is minimal in its free homotopy class. A pants decomposition is a collection of curves whose complement is a collection of three-holed spheres, often called pants, except sometimes in the United Kingdom. The length of a closed curve $\alpha$ is, by definition, the length $\ell_X(\alpha)$ of its unique geodesic minimizer on $X$. As $X$ varies in $\M_g$, so does its length. 

The geometry of a surface, in the neighborhood of a simple closed geodesic $\alpha$, is governed by its collar, or tubular neighborhood of width 
\[
w(\alpha) := \arcsinh\left(\frac{1}{\sinh\left(\frac{\ell_X(\alpha)}{2}\right)}\right).
\]
Ffor $\varepsilon>0$, the $\varepsilon$-thick part of a surface is the set of all points of injectivity radius at least $\varepsilon$. Its complement is called the $\varepsilon$-thin part. Note that all closed geodesics strictly contained in the $\varepsilon$-thick part of a surface are of length at least $2\varepsilon$.

\subsection{Useful lemmas}

Here we state a few results, mostly well known but sometimes stated differently, that will turn out to be crucial in the main proofs. 

We begin with a result about how to triangulate the thick part of a surface. We set $\varepsilon_0:= \frac{1}{2} \arcsinh(1) < \frac{1}{2}$ and consider the $\frac{\varepsilon_0}{2}$-thick part of $X$. If the corresponding thin part is empty, we set $X_0$ to be $X$. 

\begin{figure}[h]
 \centering
\begin{overpic}[width=0.4\linewidth,
tics=10]{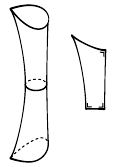}
\put(65,50){$h$}
\put(44,50){$H$}
\put(01,95){$p_i$}
\put(01,04){$q_i$}
\put(20,48){$\sigma_i$}
\put(53,76){$\varepsilon_0$}
\put(52,28){$\frac{\ell_X(\sigma_i)}{2}$}
\end{overpic}
\caption{The loop collar around $\sigma_i$}
 \label{fig:collar}
\end{figure}

If not, $X$ contains a collection of at most $3g-3$ simple  $\sigma_1,\hdots,\sigma_t$ of length at most $\varepsilon_0$. For each $\sigma_i$, we choose points $p_i,q_i$, one in the each half-collar of $\sigma$, which are the basepoints of a geodesic loops freely homotopic to $\sigma_i$ and of length $2 \varepsilon_0$ (see Figure \ref{fig:collar}). We set $X_0$ to be the surface obtained by cutting $X$ along these geodesic loops and removing the corresponding thin parts (referred to as loop-collars below).

We record two important properties of these points and associated loops: they are pairwise at distance (on $X$) at least $\varepsilon_0$ from one another. This comes from the geometry of the loop collars. To see this, consider the quadrilateral with three right angles used to construct the loop collar, as in Figure \ref{fig:collar}. One side of length $\varepsilon_0$, and the opposite side is of length $\ell_X(\sigma_i)/2$. The sides denoted $H$ and $h$ in the figure are thus both of length at least $h$ which satisfies 
\[
\cosh(h) \geq \tanh(\varepsilon_0) \coth\left( \frac{\varepsilon_0}{2}\right)
\]
and it follows that the full width of this loop collar is actually of width at least $4 \varepsilon_0$. In particular, any path that passes from one side of the loop collar to the other is also of length at least $4 \varepsilon_0$. We can now find a geometric decomposition of $X_0$. 

\begin{lemma}\label{lem:triangulation}
For $X\in \M_g$, let $\varepsilon_0$ and $X_0$ be as above. There exists a geodesic cell decomposition $\Gamma$ of $X_0$ with at most $83(g-1)$ vertices such that all edges have lengths in $[\varepsilon_0,2\varepsilon_0]$.
\end{lemma}

\begin{proof}
This follows from a standard argument which is provided for completeness. For instance, in \cite{BuserBook}, the existence of a triangulation of the thick part of a surface is given with nearly identical properties. 

For our vertex set, we first consider the points used as base points for the geodesic boundary loops of $X_0$. As noted above, they are all distance at least $\varepsilon_0$ from one another on $X$, hence on $X_0$.

We now we complete the base points of the boundary curves by considering points in $X_0$ that complete them into a maximal set of points, again, all at distance at least $\varepsilon_0$ from one another. Denote the resulting set of points $P=\{p_1,\hdots,p_N\}$. The set of balls of radius $\varepsilon_0$ around $P$ covers $X^{\varepsilon_0}$, otherwise the set is not maximal. Now observe that the balls of radius $\frac{\varepsilon_0}{2}$ around points of $P$ are embedded on $X$ and pairwise disjoint. 

Hence
\[
N\,2\pi\left(\cosh\left(\frac{\varepsilon_0}{2}\right)-1\right)
\leq \area(X)=4\pi(g-1).
\]
and therefore
\[
N\leq \frac{2(g-1)}{\cosh(\varepsilon_0/2)-1} < 83 (g-1).
\]

We now consider the Voronoi decomposition of $X_0$ associated to $P$ and its geodesic Delaunay dual. Note that, since balls of radius $\varepsilon_0$ cover, the Voronoi decomposition is a cell decomposition, hence its Delaunay dual is as well. Now since a point on the boundary of two Voronoi cells is at distance at most $\varepsilon_0$ from both points, a Delaunay edge is of length at most $2 \varepsilon_0$. This proves the claim. 
\end{proof}

The proof actually gives a triangulation, but we won't make use of that fact here. The geometric graph thus obtained is denoted by $\Gamma$. 

We now state a (very) useful lemma from the world of circle graphs. 

\begin{lemma}\label{lem:circle}
Let $n \geq 2 $ pairs of distinct points lie on a circle. Two pairs are said to alternate if their endpoints alternate on the circle. There is a set of at least
\[
 \sqrt{\frac{n}{\lceil\log_2(n)\rceil}}
\]
pairs which either all mutually alternate or never alternate with one another.\end{lemma}

\begin{proof}
This follows from a bound on the chromatic number of circle graphs. More precisely, form a graph $G$ with $n$ vertices, one for each pair of points and join them by an edge if they alternate. Consider the complement graph $\overline{G}$ where you join pairs if they do not alternate. In terms of these graphs, the lemma is a lower bound on the size of their maximal complete subgraphs, that is their clique number. 

Now in \cite{Cerny2007} it is shown that 
\[
\chi\left( G\right) \leq  \omega\left(G \right) \lceil \log_2(n)\rceil
\]
where $\chi\left(G\right)$ is the graph's chromatic number and $\omega\left( G \right) $ is its clique number.  Since a proper coloring of $G$ partitions $G$ into $\chi(G)$ independent sets, and that $\omega\left( \overline{G}\right)$ is the independence number of $G$, if follows that
\[
n \leq \omega(\overline{G}) \chi(G) \leq  \omega(\overline{G}) \omega\left(G \right) \lceil \log_2(n)\rceil
\]
and hence the maximum of the clique numbers of $\overline{G}$ and $G$ satisfies the desired bound.
\end{proof}
Although the estimate used here is explicitly stated in \cite{Cerny2007}, the same order of growth is already established in \cite{Kostochka1988} where examples are exhibited showing that you cannot hope for a better order of growth. 

The way this lemma applies in our setting will become clear, but here is a quick spoiler. We will have disjoint arcs with endpoints on the boundary curve of a surface with topology. In order to construct curves with these arcs, we will have to know whether their endpoints alternate or not. 

\section{Breaking the surface into spheres}

In this part, we establish the most difficult part of the bound: the existence of a set of $g$ curves of total length at most $\approx g^{3/2} \sqrt{\log(g)}$ such that their complement on $X$ is a sphere with $2g$ boundary curves (Theorem \ref{thm:cutsystemintro} from the introduction). 

\begin{remark}\label{rem:combinatorial}
While the goal is to prove this for hyperbolic surfaces, it is really a proof that applies to a graph embedding. 

More specifically, we consider $X$ and the thick part $X_0$ of $X$. We have shown the existence of an embedded geometric graph $\Gamma$ that captures the geometry of $X_0$. The entire surface $X_0$ is very close to $\Gamma$, so we can think of it as being a skeleton of the surface. The geometry of $X_0$ is really governed by the geometry of $\Gamma$. 

In particular, all curves on $X_0$ are freely homotopic to closed paths on $\Gamma$. This gives two ways of measuring the length of a curve: by either the  length in of a minimal representative of the curve on graph $\Gamma$ or by the length of the corresponding closed geodesic on $X$. 

For curves $\alpha$ constructed below, we denote by $\ell_\Gamma(\alpha)$ the combinatorial length of the particular closed path on $\Gamma$ supplied by the construction, where edges have length $1$ and are counted with multiplicity. Since every
edge of $\Gamma$ has length at most $2\varepsilon_0$, we have
\[
\ell_X(\alpha)\leq 2\varepsilon_0\,\ell_\Gamma(\alpha).
\]
In particular, the methods used here apply to a more general combinatorial setting, this whenever you have a cellular graph embedding with a bound on the topology of the surface and the length of the graph, the method gives you a bounded pants decomposition of controlled length. 
\end{remark}

We now proceed to the proof. Consider $X_0$ and $\Gamma$ constructed in the previous section. If $X_0$ is planar, there is nothing to prove, so we suppose that the genus $g(X_0)$ of $X_0$ is at least $1$. 

While $\Gamma$ might not be connected, we work simultaneously on all of its connected components. In first instance, we proceed as if $\Gamma$ was connected and in Remark \ref{rem:disconnect}, we explain how to modify the argument to the case of multiple connected components. 

We construct a type of spanning tree $\tau$ with the set of leaves that contain the boundary loops. To do this, as a graph, we replace the boundary loops by a single edge we call a twig. (This is the same thing as adding a vertex in the middle of the loop and collapsing the two edges.) We consider an actual spanning tree of the graph thus obtained, on which the twigs are necessarily leaves, and then we put back the boundary loops to obtain $\tau$. The boundary of a small neighborhood of $\tau$ is a simple closed curve we denote $\beta_0$, which may or may not be essential, or have a geodesic representative on $X_0$, depending on whether $X_0$ has boundary or not. 

We make an important observation. By considering a minimal length representative on $X$, we have $\ell_X(\beta_0) \leq 2 \ell_X(\tau)$. We can think of $\beta_0$ lying in a very small neighborhood of $\tau$. We consider the surface $X_1$ obtained by cutting $X_0$ along $\beta_0$: it separates a planar surface $Y_0$ (which in our world is very thin) from a surface $X_1$, also of genus $g_0$ and with a single boundary curve $\beta_0$. And in fact, all our curves, and later planar subsurfaces, will lie in neighborhoods of $\Gamma$. For good measure, we set $g_1= g(X_1) = g_0$.

We think of $\tau$ as a subgraph of $\Gamma$ which uses a certain number, but not all of the arcs of $\Gamma$. Let $\arc$ be set of arcs in $\Gamma \setminus \tau$. Each of these arcs also live on $X_1$, albeit in a slightly truncated version, with endpoints on $\beta_0$. They may or may not be pairwise homotopic, relative to endpoints sliding on $\beta_0$. (In fact, because $\Gamma$ was originally a triangulation, their homotopy classes provide a full hexagon decomposition of $X_1$, see for instance \cite{GultepeParlier2026}, but we won't use that here.) 

We now introduce a standard trick, used for instance in \cite{APP2011,DisarloParlier2018,DisarloParlier2019}. We provide a sketch proof for completeness.

\begin{lemma}\label{lem:arcs}
The set of arcs $\arc$ contains a subset $\arc_1$, consisting of $2g_1$ arcs, such that $X_1 \setminus \arc_1$ is a disk.
\end{lemma}

\begin{figure}[htbp]
\centering
\begin{subfigure}[b]{0.45\linewidth}
  \centering
  \begin{overpic}[width=\linewidth,tics=10]{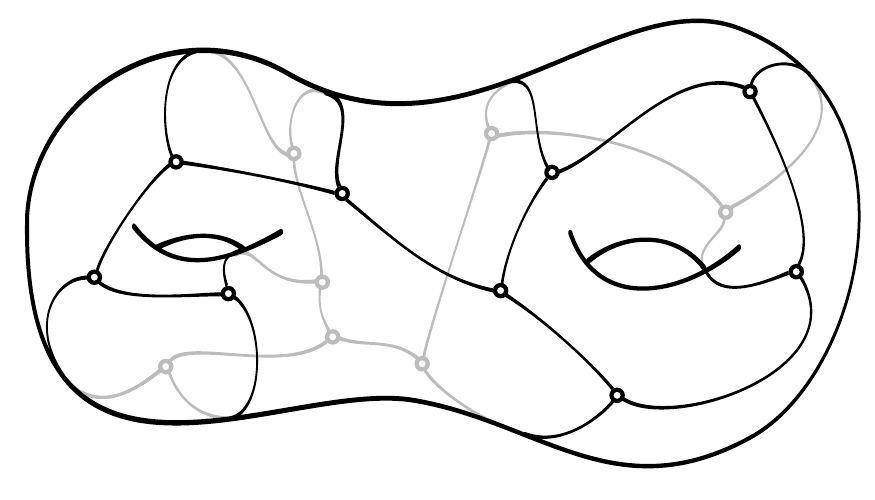}
     \end{overpic}
  \caption{The cell decomposition $\Gamma$}
  \label{fig:a}
\end{subfigure}
\hfill
\begin{subfigure}[b]{0.45\linewidth}
  \centering
  \begin{overpic}[width=\linewidth,tics=10]{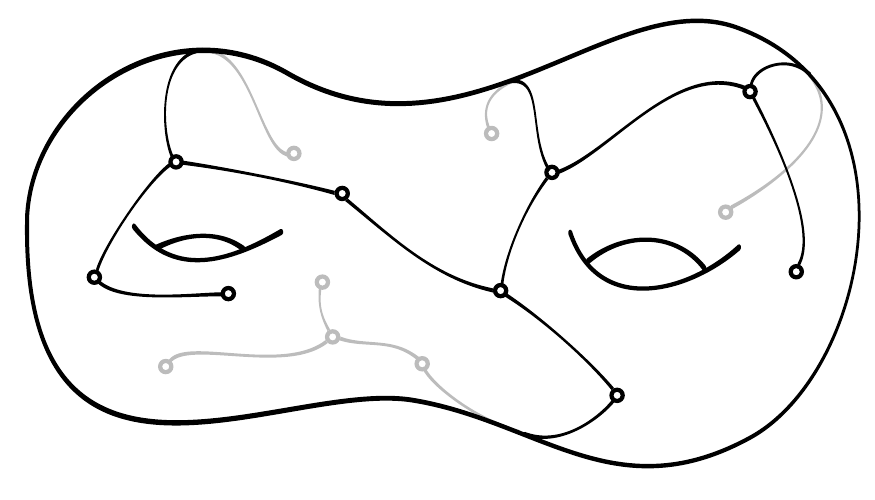}
  \end{overpic}
  \caption{A spanning tree of $\Gamma$}
  \label{fig:b}
\end{subfigure}

\vspace{1em}

\begin{subfigure}[b]{0.45\linewidth}
  \centering
  \begin{overpic}[width=\linewidth,tics=10]{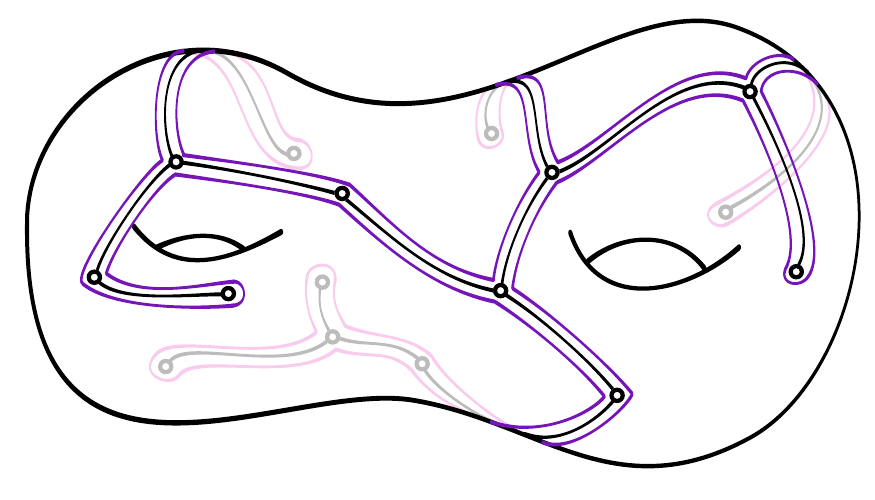}
  \end{overpic}
  \caption{The curve $\beta_0$}
  \label{fig:c}
\end{subfigure}
\hfill
\begin{subfigure}[b]{0.45\linewidth}
  \centering
  \begin{overpic}[width=\linewidth,tics=10]{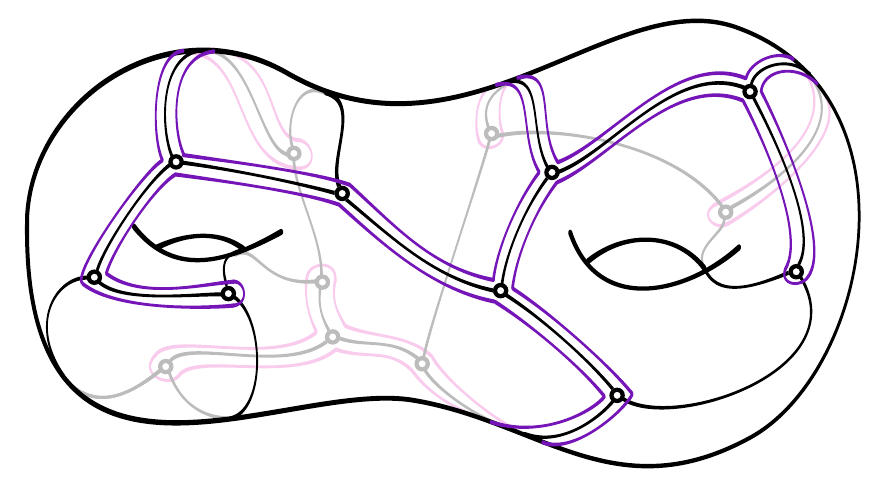}
  \end{overpic}
  \caption{Putting back the arcs}
  \label{fig:d}
\end{subfigure}

\vspace{1em}

\begin{subfigure}[b]{0.45\linewidth}
  \centering
  \begin{overpic}[width=\linewidth,tics=10]{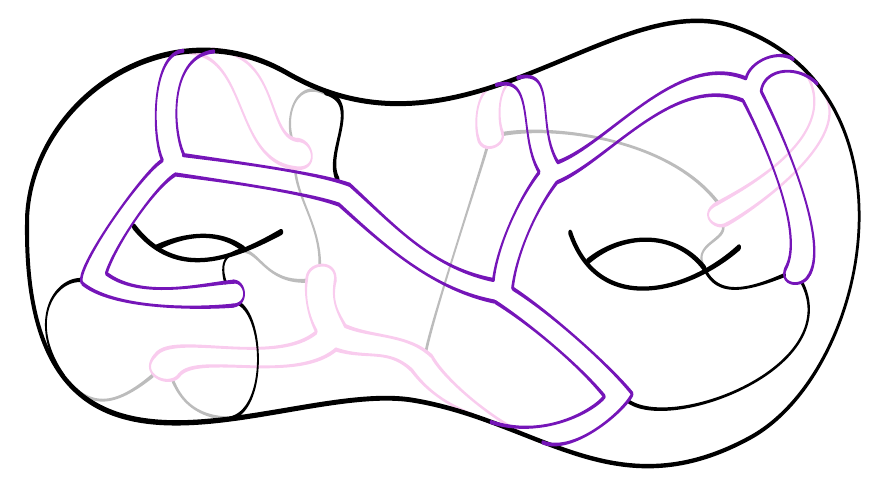}
  \end{overpic}
  \caption{Forgetting the spanning tree}
  \label{fig:e}
\end{subfigure}
\hfill
\begin{subfigure}[b]{0.45\linewidth}
  \centering
  \begin{overpic}[width=\linewidth,tics=10]{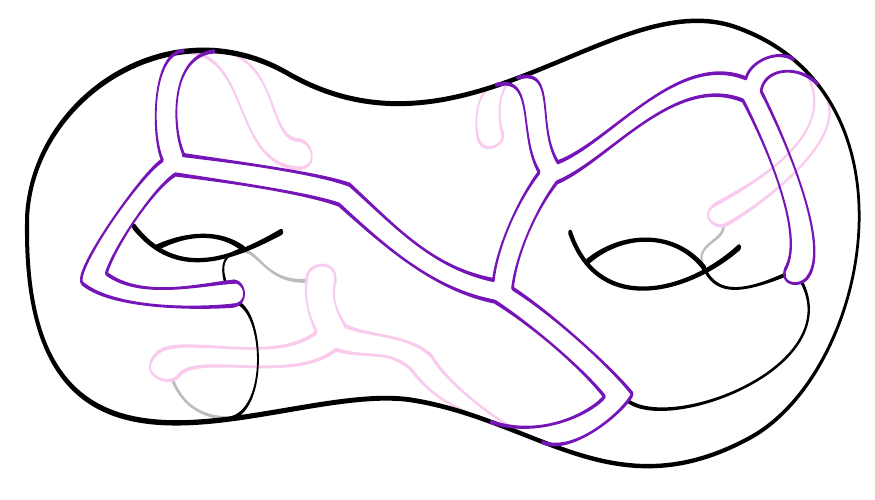}
  \end{overpic}
  \caption{Choosing a set of active arcs}
  \label{fig:f}
\end{subfigure}
\caption{An example of the construction of the active arcs $\arc_1$ in genus $2$. In this example, $\tau$ does not contain any loops.}
\label{fig:thefullmonty}
\end{figure}

\begin{proof}
Since $X_1 \setminus \arc$ is a collection of simply connected pieces, it is of genus $0$. Hence, there must exist an arc $a_1\in \arc$ such that $X_1 \setminus a_1$ is of genus $g_0 - 1$. This results (topologically) in a surface with two boundary components. Now there must be an arc $a_2$ which links these two boundary components. We cut along this arc we are now back to a single boundary component and a connected region of genus $g_0 -1$. This procedure can be iterated until we reach a set of arcs $\arc_1$ such that their complementary region is a surface of genus $0$ with a single boundary component, generally referred to as a disk. 
\end{proof}
An example of the entire construction of $\arc_1$ in genus $2$, starting from a geometric cell decomposition, is illustrated in Figure \ref{fig:thefullmonty}.

The set of arcs $\arc_1$ will be referred to as our initial set of {\it active arcs}, and as we use them to construct curves and subsurfaces, they will become inactive. Note that cutting along any subset of them never disconnects the surface, otherwise it would be impossible to have a single disk by cutting along all of them. 

Also note that each arc has distinct endpoints on $\beta_0$, hence define a pair of points on a circle as in Lemma \ref{lem:circle}. And in fact, we will inactivate them in batches, inductively, where the size of the inactivated batch will come from this lemma.

\subsection{The inductive step.}\label{ss:induct}

We begin at $k=1$. 

Set $n_k:=| \arc_k |$. By Lemma \ref{lem:circle}, the set $\arc_k$ contains at least $ \sqrt{\frac{n_k}{\lceil\log_2(n_k)\rceil}}$ arcs whose endpoints either mutually alternate or mutually do not alternate. We set $m_k$ to be the largest even integer such that there are $m_k$ arcs whose endpoints either mutually alternate or mutually do not alternate. It satisfies
\[
m_k \geq \sqrt{\frac{n_k}{\lceil\log_2(n_k)\rceil}} - 2.
\]
We begin with the case where they do not alternate.

\noindent{\it Non-alternating arcs.}

We consider the $m_k$ arcs, which we denote by $a_1 \hdots,a_{m_k}$, that pairwise do not alternate. (Formally, we should write something like $a^{[k]}_1 \hdots,a^{[k]}_{m_k}$ because these are the arcs we use at stage $k$, and we will get rid of them at the end of the iteration, but for simplicity we omit the superscript.) We consider the subsurface $N_k$ obtained by considering a small neighborhood of $\beta_{k-1}$ together with the arcs. It is planar and has $m_k+2$ boundary curves. One of them is $\beta_{k-1}$. 

\begin{figure}[h]
 \centering
\begin{overpic}[width=0.4\linewidth,
tics=10]{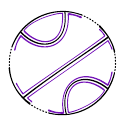}
\put(33,71){$\alpha^{[k]}_{i}$}
\put(22,54){$\alpha^{[k]}_{j}$}
\end{overpic}
\caption{Constructing the curves $\alpha^{[k]}_i$}
 \label{fig:nonalt}
\end{figure}

The remaining $m_k+1$ boundary curves we denote by
\[
\alpha^{[k]}_1,\hdots, \alpha^{[k]}_{m_k+1}.
\]
Note that some may be trivial but certainly not all, see Figure \ref{fig:nonalt}). We add $N_k$ to our collection of planar subsurfaces and its boundary curves to our collection of curves. 

Observe that the complementary region on $X_k$ to $N_k$, denoted $R_k$, has $m_k+1$ boundary curves and has genus $g_{k} - m_k$ because each arc reduces the genus by exactly $1$. It is also connected and, in fact, the remaining arcs connect these boundary curves. (If not, it would be impossible to cut along a collection of remaining arcs to obtain a single disk.) 

We can now chose $m_k$ of the remaining arcs, $a_{m_k+1},\hdots,a_{2m_k} \in \arc_k$, so that cutting along them reduces the number of boundary components to $1$. We call them {\it connecting arcs}. We consider the planar surface $M_k$ obtained by taking a neighborhood of $\partial R_k$ together with the arcs $a_{m_k+1},\hdots,a_{2m_k}$. 

The complement of $M_k$ in $R_k$ is a connected surface, which we denote by $X_{k+1}$, with one boundary curve $\beta_k$ and genus $g_{k+1}=g(X_{k+1})=g_k-m_k.$

We remove all of the arcs we used from the set of active arcs to obtain a new set of active arcs $\arc_{k+1}$. 

\noindent{\it Alternating arcs.}

We consider the $m_k$ arcs, $a_1 \hdots,a_{m_k}$, that pairwise alternate (again we've used the simpler notation which omits superscripts).

We consider the subsurface $S_k$ obtained by considering a small neighborhood of $\beta_{k-1}$ together with the arcs. This time it is not planar, and in fact has genus $m_k/2$, and, due to the parity of $m_k$, two boundary curves. The genus can be computed via a straightforward Euler characteristic computation. 

One of the boundary curves is $\beta_{k-1}$ and the other we denote by $\beta_k$. We will describe the homotopy class more precisely in what follows.

We give $\beta_{k-1}$ an orientation and choose an endpoint of $a_1$. Up to reordering if necessary, we can suppose that the next endpoint we meet along $\beta_{k-1}$ belongs to  $a_2$ and, more generally, for $i=1,\hdots,m_k-1$, after the first endpoint of $a_i$ the next one along $\beta_{k-1}$ is  $a_{i+1}$. Since the arcs all alternate, after $a_{m_k}$, the next endpoint belongs to $a_1$ and the pattern repeats. This also gives a way to give the arcs orientation: we orient the arcs to go between their first and second endpoints.

\begin{figure}[h]
 \centering
\begin{overpic}[width=0.4\linewidth,
tics=10]{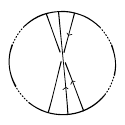}
\put(55,16){$a_1$}
\put(63,25){$a_{2}$}
\put(59,78){$a_{m_k}$}
\put(59,03){$b_1$}
\put(50,94){$b_{m_k}$}
\put(38,93){$\bar{b}_1$}
\end{overpic}
\caption{The arcs $a_i$, $b_i$ and $\bar{b}_i$}
 \label{fig:alt01}
\end{figure}

For $i=1,\hdots,m_k-1$, we set $b_i$ to be the oriented arc of $\beta_{k-1}$ between the first endpoints of $a_i$ and $a_{i+1}$ and $b_{m_k}$ to be between the first endpoint of $a_{m_k}$ and the second of $a_{1}$. For $i=1,\hdots,m_k-1$, we set $\bar{b}_i$ to be the oriented arc between the second endpoints of $a_i$ and $a_{i+1}$ and $\bar{b}_{m_k}$ to be the arc of $\beta_{k-1}$ between the second endpoint of $a_{m_k}$ and the first of $a_{1}$. See Figure \ref{fig:alt01}.

\begin{figure}[h]
 \centering
\begin{overpic}[width=0.4\linewidth,
tics=10]{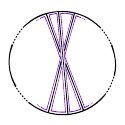}
\put(63,34){$\beta_{k}$}
\put(76,85){$\beta_{k-1}$}
\end{overpic}
\caption{Constructing $\beta_{k}$ (in color)}
 \label{fig:betak}
\end{figure}

The curve $\beta_k$ is in the homotopy class of the cyclic concatenation
\[
\left(
\prod_{j=1}^{m_k/2}
a_{2j-1}*\bar b_{2j-1}*a_{2j}^{-1}*b_{2j}
\right)
*
\left(
\prod_{j=1}^{m_k/2}
a_{2j-1}^{-1}*b_{2j-1}*a_{2j}*\bar b_{2j}
\right),
\]
where the products denote concatenation in increasing order of $j$.
The final endpoint is the starting point, so this concatenation
is closed. See Figure \ref{fig:betak}.

We will now construct a collection of disjoint curves that cut $S_k$ into a sphere with $m_k+2$ holes. Again, the parity of $m_k$ is essential here. 

For each even $i=2,4,\hdots,m_k$, we set $\alpha^{[k]}_i$ to be the curve in the free homotopy class of the cyclic concatenation of $a_i$, $\bar{b}_{i-1}^{-1}$, $a_{i-1}^{-1}$ and $b_{i-1}$. These $\frac{m_k}{2}$ curves are disjoint as they do not share any arcs, and cutting along all the curves $\alpha^{[k]}_2, \hdots, \alpha^{[k]}_{m_k}$ results in a connected surface. As the genus of $S_k$ is $\frac{m_k}{2}$, this resulting connected surface $Y_k$ is necessarily a sphere with $2+ m_k$ boundary curves. We add the resulting surface $Y_k$ to our collection of planar surfaces, and the curves we cut along to our collection of curves. 

The closure of the complement of $S_k$ in $X_k$ is a connected surface $X_{k+1}$ with boundary $\beta_k$ and genus
\[
g_{k+1}=g_k-\frac{m_k}{2}.
\]

We remove all of the arcs we used from the set of active arcs to obtain a new set of active arcs $\arc_{k+1}$.

In both the alternating and non-alternating cases, $X_{k+1}$ is obtained by cutting $X_k$ along the arcs just used. Consequently, the remaining active arcs $\arc_{k+1}$ cut $X_{k+1}$ into a disk, so the construction can be iterated.

\subsection{Controlling the length of curves and the number of batches.}

\subsubsection{Bounding the number of stages}
Recall that $n_k:=| \arc_k |$ and that $n_1 = 2 g_0 \leq 2 g$. At each stage, we reduced the size of $n_k$ by at least $m_k$ where 
\[
m_k \geq \sqrt{\frac{n_k}{\lceil\log_2(n_k)\rceil}} - 2,
\]
hence
\[n_{k+1} \leq n_k - \sqrt{\frac{n_k}{\lceil\log_2(n_k)\rceil}} + 2.
\]
Since \(m_k\) is positive and even at each stage, we have
\[
n_k-n_{k+1}\geq m_k
\geq \max\left\{2,
\sqrt{\frac{n_k}{\lceil\log_2 n_1\rceil}}-2\right\}
\geq \frac{\sqrt{n_k}}{2\sqrt{\lceil\log_2 n_1\rceil}}.
\]
Consequently,
\[
\sqrt{n_k}-\sqrt{n_{k+1}}
=\frac{n_k-n_{k+1}}{\sqrt{n_k}+\sqrt{n_{k+1}}}
\geq \frac{1}{4\sqrt{\lceil\log_2 n_1\rceil}}.
\]
Summing over all stages, the left-hand side telescopes to
\(\sqrt{n_1}\), since the terminal value is \(0\).
Thus the number of stages is at most
\[
4\sqrt{n_1\lceil\log_2 n_1\rceil}
\leq 4\sqrt{2n_1\log_2 n_1}.
\]
\subsubsection{The lengths of a cut system}
We first make an observation about the length of the curves $\beta_k$ that we have constructed. 

The curve $\beta_0$ traverses each edge of the spanning tree twice and each boundary loop once. The spanning tree has $N-1$ edges, and there are at most $6g-6$ boundary loops so
\begin{equation}\label{eq:beta0}
\ell_\Gamma (\beta_0) <172 (g-1).
\end{equation}

\begin{lemma}\label{lem:betacurves}
For $k\geq 1$, the curves $\beta_k$ satisfy the following bound:
\[
\ell_\Gamma (\beta_k) \leq \ell_\Gamma(\beta_0) + 4 |\arc_1 |.
\]
\end{lemma}
\begin{proof}
This follows from the construction of the curves $\beta_k$. In either the alternating or non-alternating case, they always satisfy the bound
\[
\ell_\Gamma (\beta_k) \leq \ell_\Gamma (\beta_{k-1}) + 4 m_k.
\]
In the alternating case, this comes from the fact that $\beta_k$ is constructed using all arcs of $\beta_{k-1}$ once and twice each alternating arc (so in fact satisfies the stronger inequality with a $2$ instead of a $4$).

Now in the non-alternating case, we first constructed the curves $\alpha^{[k]}_1,\hdots,\alpha^{[k]}_{m_k+1}$ using $\beta_{k-1}$ and each of the $m_k$ arcs twice. We then used $m_k$ additional connecting arcs, again twice, to construct $\beta_k$. Hence
\[
\ell_\Gamma (\beta_k) \leq \sum_{i=1}^{m_k+1}\ell_\Gamma \left( \alpha^{[k]}_i\right) +  2 m_k \leq \ell_\Gamma (\beta_{k-1}) + 4 m_k.
\]
Furthermore, at the end of each stage, we remove the used arcs from $\arc_{k}$ to obtain $\arc_{k+1}$. In the non-alternating case we have used $2m_k$ arcs, and in the alternating case only $m_k$. So in either case we have
\[
|\arc_{k+1}| \leq |\arc_{k}| - m_k
\]
and thus
\[
m_k \leq  |\arc_{k}| -|\arc_{k+1}|.
\]
By combining the two inequalities we obtain 
\[
\ell_\Gamma (\beta_k)  \leq \ell_\Gamma (\beta_{k-1}) + 4 \left( |\arc_{k}| -|\arc_{k+1}| \right) 
\]
which implies, using telescoping sums, that 
\begin{eqnarray*}
\ell_\Gamma (\beta_k) & \leq &\ell_\Gamma (\beta_{k-2}) + 4 \left( |\arc_{k-1}| -|\arc_{k+1}|\right) \\
& \leq &\ell_\Gamma (\beta_{0}) + 4 \left( |\arc_{1}| -|\arc_{k+1}|\right)\\
& \leq &\ell_\Gamma (\beta_{0}) + 4 |\arc_{1}|
\end{eqnarray*}
as required.
\end{proof}

We used the following inequality on the sum of the lengths of the curves used in the non-alternating case, which we record here:

\begin{lemma}\label{lem:nonaltcurves}
The curves $\alpha^{[k]}_1,\hdots,\alpha^{[k]}_{m_k+1}$ constructed in the non-alternating case satisfy
\[
\sum_{i=1}^{m_k+1}\ell_\Gamma\left(\alpha^{[k]}_i\right)
\leq \ell_\Gamma(\beta_{k-1})+2m_k
\leq \ell_\Gamma(\beta_0)+4|\arc_1|.
\]

\end{lemma}
\begin{proof}
This is a direct consequence of the construction of the curves $\beta_{k-1}$ and the previous telescoping proof. 
\end{proof}

In the alternating case, we can also bound the lengths of the curves $\alpha^{[k]}_2,\hdots,\alpha^{[k]}_{m_k}$ (there are $m_k/2$ of them) using the same estimates. 

\begin{lemma}\label{lem:altcurves}
The curves $\alpha^{[k]}_2,\hdots,\alpha^{[k]}_{m_k}$ satisfy
\[
\sum_{i=1}^{m_k/2}\ell_\Gamma \left( \alpha^{[k]}_{2i}\right) \leq \ell_\Gamma (\beta_{k-1}) + m_k \leq  \ell_\Gamma (\beta_{0}) + 4 |\arc_{1}|.
\]
\end{lemma}

In total, if the process ends after $s$ stages, the total length of a multicurve such that cutting $X_0$ along this multicurve is of genus $0$ is at most $s$ times the max lengths of the curves $\beta_k$s, and $s$ times the sums of the lengths of the curves $\alpha^{[k]}_i$s. Using the previous lemmas, this set of curves is of total length at most
\[
s \left( 2 \ell_\Gamma (\beta_{0}) + 16 g \right) 
\] 
and thus the total length is at most $360 gs$. From this we can deduce the following:
\begin{lemma}
The surface $X_0$ admits a cut system of total length at most $360 gs $.
\end{lemma}
\begin{proof}
The curves we have constructed can all be realized disjointly, as simple curves, and with lengths strictly below the claimed bound. While some of them might be homotopically trivial or separating, the set certainly contains a curve that reduces the genus by exactly $1$. By iterating this argument, we get $g_0$ curves which, by cutting along them, reduces the genus to $0$ and hence form a cut system. 
\end{proof}

\begin{remark}\label{rem:disconnect}
If $\Gamma$ is disconnected, we apply the construction simultaneously to each of its components. The length estimates add over the components and the number of stages is exactly the maximum of the number of stages needed for each component. As such, the connected case already accounts for the worst case scenario, and the same bounds hold without an assumption on the connectedness of $\Gamma$. 
\end{remark}

We now promote the cut system to one of $X$ by taking into consideration the curves we cut $X$ along to obtain $X_0$.

\begin{theorem}\label{thm:cutsystem}
Every closed hyperbolic surface $X$ of genus $g\geq 2$ admits a cut system of total length at most
\[
3000\,g^{3/2}\sqrt{\log g},
\]
where $\log$ is the natural logarithm. 
\end{theorem}

\begin{proof}
Take the cut system of $X_0$, and add to them the core geodesics of the removed collars. The resulting multicurve has planar complement, so, using the same argument as previously, it contains a cut system $\alpha_1, \hdots,\alpha_g$ of $X$. 

At most $g$ of its curves are core curves of the thin collars, and each has length at most $2\varepsilon_0$. Since every edge of $\Gamma$ also has length at most $2\varepsilon_0=\arcsinh(1)$, the preceding combinatorial estimates give
\[
\ell_X( \alpha_1 \cup  \hdots \cup \alpha_g)
\leq \arcsinh(1) \left( (360g-344)s+g \right)
\]
where $s$ is the number of stages in the combinatorial construction. Now using the bound on $s$, and since $360g-344\leq360\sqrt{g(g-1)}$, it follows that
\[
\ell_X( \alpha_1 \cup  \hdots \cup \alpha_g)
\leq \arcsinh(1)\left(
1440g\sqrt{2(g-1)\lceil\log_2(2g)\rceil}+g
\right).
\]
Now using elementary inequalities, we obtain
\[
\ell_X( \alpha_1 \cup  \hdots \cup \alpha_g) < 3000\,g^{3/2}\sqrt{\log g}
\]
as required.
\end{proof}

\subsection{Pants decompositions}

We use the preceding results to deduce an upper bound on the sums of lengths of a minimal length pants decomposition. 

We first prove a technical lemma which is essentially already contained in \cite{BPS2012}. 

\begin{lemma}\label{lem:binary}
Let $Y$ be a planar component obtained by cutting $X$ along disjoint essential curves, with $b\geq3$ boundary components. Let $\Gamma_Y \subset Y$ be a connected embedded graph containing $\partial Y$. Then $Y$ admits a pants decomposition $P$ such that
\[
\sum_{\alpha\in P}\ell_X(\alpha)
\leq4\lceil\log_2 b\rceil\,\ell_X(\Gamma_Y).
\]
\end{lemma}

\begin{proof}
Choose a tree joining the boundary components after collapsing each boundary circle to a vertex. Following its boundary orders the boundary components. Recursively bisect this order, as in \cite[proof of Proposition~6.8]{BPS2012}, and take curves surrounding the resulting blocks.

At each level, each tree edge is traversed at most four times and each original boundary edge at most three times. Thus the total length at that level is at most $4\ell_X(\Gamma_Y)$. There are at most $\lceil\log_2 b\rceil$ levels. Discarding peripheral and repeated curves and taking geodesic representatives proves the estimate.
\end{proof}

We can now prove the main result, Theorem \ref{thm:pantsintro} from the introduction. 

\begin{theorem}\label{thm:pants}
Every closed hyperbolic surface $X$ of genus $g\geq2$ admits a pants decomposition of total length at most
\[
80000\,g^{3/2}(\log g)^{3/2}.
\]
\end{theorem}

\begin{proof}
Let $M$ be the full multicurve constructed on $X_0$, after discarding trivial curves, boundary-parallel curves and redundant parallel copies. We set $L=\ell_X(M)$.
The length calculation right before Theorem \ref{thm:cutsystem} bounds the full collection before choosing a subset that formed a cut system. In particular,
\[
L\leq \arcsinh(1)(360g-344)s
\leq 3000g^{3/2}\sqrt{\log g}.
\]
Moreover, the vertex and edge bounds for $\Gamma$ give
\[
\ell_X(\Gamma)\leq255(g-1).
\]
Now if we cut $X_0$ along $M$, we obtain planar surfaces $Y_i$, each with 
$b_i \leq 2 g$ boundary components.

Cutting $\Gamma$ and adjoining the new boundary curves gives a connected cellular graph $\Gamma_i$ on each $Y_i$, containing its boundary. We have
\[
\sum_i  \ell_X(\Gamma_i) \leq 255g+2L.
\]
We can now apply Lemma~\ref{lem:binary} to each $Y_i$. Together with $M$ and the core geodesics of the removed collars, the resulting curves give a pants decomposition
$P$ of $X$ after taking geodesic representatives. The core geodesics contribute less than $3g$, and hence
\[
\ell_X(P)
\leq L+3g+ 4\lceil\log_2(2g)\rceil(255g+2L).
\]
Now, using $g\geq 2$ and elementary estimates, we obtain
\[ \ell_X(P) <80000\,g^{3/2}(\log g)^{3/2}
\]
as claimed. 
\end{proof}

{\it Address:}\\
Department of Mathematics, University of Fribourg, Switzerland\\

{\it Email:}\\
hugo.parlier@unifr.ch\\

\end{document}